\documentclass[12pt,a4paper,english,intoc,bibliography=totoc,index=totoc,BCOR10mm,captions=tableheading,titlepage,fleqn]{scrbook}
\usepackage{lmodern}

\usepackage[T1]{fontenc}
\usepackage[latin9]{inputenc}
\usepackage{babel}
\usepackage{amsmath}
\usepackage{amssymb}
\usepackage[unicode=true,
 bookmarks=true,bookmarksnumbered=true,bookmarksopen=true,bookmarksopenlevel=2,
 breaklinks=false,pdfborder={0 0 0},backref=false,colorlinks=false]
 {hyperref}
\hypersetup{pdftitle={Abstract},
 pdfauthor={Uwe Stöhr},
 pdfpagelayout=OneColumn, pdfnewwindow=true, pdfstartview=XYZ, plainpages=false}
\usepackage{breakurl}

\makeatletter

\usepackage{calc}

\makeatother

\begin{document}
\noindent \begin{center}
\textbf{\Large Transcendence bases, well-orderings of the reals and
the axiom of choice}
\par\end{center}{\Large \par}

\noindent \begin{center}
{\large Haim Horowitz and Saharon Shelah}
\par\end{center}{\large \par}

\noindent \begin{center}
\textbf{\small Abstract}
\par\end{center}{\small \par}

\noindent \begin{center}
{\small We prove that $ZF+DC+"$there exists a transcendence basis
for the reals$"+"$there is no well-ordering of the reals'' is consistent
relative to $ZFC$. This answers a question of Larson and Zapletal.}%
\footnote{{\small Date: January 27, 2019}{\small \par}

2010 Mathematics Subject Classification: 03E25, 03E35, 03E40, 12F20

Keywords: transcendence basis, well-ordering, axiom of choice, forcing,
amalgamation

Publication 1093 of the second author

Partially supported by European Research Council grant 338821.%
}
\par\end{center}{\small \par}

\textbf{\large Introduction}{\large \par}

It's well-known that the axiom of choice has far-reaching consequences
for the structure of the real line. Among them, to name a few, are
the existence of non-measurable sets of reals, nonprincipal ultrafilters
on $\omega$, paradoxical decompositions of the unit sphere, mad families
and more. As the aforementioned statements are consistently false
over $ZF+DC$, it's natural to study the possible implications between
them in the absence of choice. This direction of study has gained
considerable interest in recent years, with many consistency results
showing mostly the independence over $ZF+DC$ between various properties
of the real line implied by the axiom of choice. We mention several
such examples:

\textbf{Theorem ({[}Sh:218{]}): }It's consistent relative to an inaccessible
cardinal that $ZF+DC$ holds, all set of reals are Lebesgue measurable
and there is a set of reals without the Baire property.

\textbf{Theorem ({[}HwSh:1113{]}): }It's consistent relative to an
inaccessible cardinal that $ZF+DC$ holds, all sets of reals are Lebesgue
measurable and there is a mad family.

\textbf{Theorem ({[}LaZa1{]}): }It's consistent relative to a proper
class of Woodin cardinals that there exists a mad family and there
are no $\omega_1$ sequences of reals, nonatomic measures on $\omega$
and total selectors for $E_0$.

Our current paper will focus on two consequences of the axiom of choice
for the real line, namely the existence of a transcendence basis for
the reals and the existence of a well-ordering of the reals. The following
question was asked by Larson and Zapletal in their forthcoming book:

\textbf{Question ({[}LaZa2{]}): }Does the existence of a transcendence
basis for the reals imply the existence of a well-ordering of the
reals?

We shall prove that the answer is negative, namely:

\textbf{Main result: }$ZF+DC+"$there exists a transcendence basis
for the reals$"+"$there is no well-ordering of the reals$"$ is consistent
relative to $ZFC$.

The proof strategy will be similar to that of {[}Sh:218{]} and {[}HwSh:1113{]}
(though no inaccessible cardinals will be used in the current proof).
Our forcing $\mathbb P$ will consist of conditions $p=(u_p,\mathbb{Q}_p,\underset{\sim}{R_p})$
where $\mathbb{Q}_p$ is a ccc forcing from some fixed $H(\lambda)$
that forces $MA_{\aleph_1}$ and $\underset{\sim}{R_p}$ is a set
of $\mathbb{Q}_p$-names of reals that's forced by $\mathbb{Q}_p$
to be a transcendence basis for the reals. The order will be defined
naturally. The sets of the form $\underset{\sim}{R_p}$ will approximate
a transcendence basis in the final model, while the forcing notions
$\mathbb{Q}_p$ will help us to prove the non-existence of a well-ordering
of the reals using a standard amalgamation argument. The fact that
each $\mathbb{Q}_p$ forces $MA_{\aleph_1}$ will guarantee that the
relevant amalgamation will be ccc.

\textbf{Acknowledgement:} We would like to thank Jindra Zapletal for
informing us about a gap in a previous version of this paper.

The rest of the paper will be devoted to the proof of the main result
mentioned above. We shall assume basic familiarity with amalgamation
of forcing notions (see, e.g., {[}HwSh:1090{]}).

\textbf{\large Proof of the main result}{\large \par}

\textbf{Hypothesis 1: }Throughout the paper, we fix infinite regular
cardinals $\lambda$ and $\kappa$ and an infinite cardinal $\mu$
such that $\mu=\mu^{\aleph_1}<\lambda$, $\kappa=\mu^+$ or $\aleph_2 \leq cf(\kappa) \leq \kappa \leq \lambda$
and $(\forall \alpha< \kappa)([\alpha]^{\aleph_1}<\kappa)$.

\textbf{Definition 2: }We define the forcing notion $\mathbb P$ as
follows:

A. $p\in \mathbb P$ iff $p=(u,\mathbb Q,\underset{\sim}{R})=(u_p,\mathbb{Q}_p,\underset{\sim}{R}_p)$
where:

a. $u\in [\lambda]^{<\kappa}$.

b. $\mathbb Q \in H(\lambda)$ is a ccc forcing such that $u$ is
its underlying set of elements.

c. $\Vdash_{\mathbb Q} MA_{\aleph_1}$.

d. $\underset{\sim}{R}$ is a set of $\mathbb Q$-names of reals that
is forced by $\mathbb Q$ to be a transcendence basis of the reals.

B. $p\leq_{\mathbb P} q$ iff

a. $u_p \subseteq u_q$.

b. $\mathbb{Q}_p \lessdot \mathbb{Q}_q$.

c. $\underset{\sim}{R}_p \subseteq \underset{\sim}{R}_q$.

\textbf{Definition 3: }We define the following $\mathbb P$ names:

a. $\underset{\sim}{ \mathbb Q}=\cup \{ \mathbb{Q}_p : p\in \underset{\sim}{G}_{\mathbb P}\}$.

b. $\underset{\sim}{R}=\cup \{ \underset{\sim}{R}_p : p \in \underset{\sim}{G}_{\mathbb P}\}$.

\textbf{Claim 4: }a. $\mathbb P$ is a forcing notion of cardinality
$\lambda^{<\kappa}$, preserving cardinals and cofinalities of cardinals
$\leq \kappa$ and $> \lambda^{<\kappa}$.

b. If $\delta<\kappa$ is a limit ordinal and $\bar p=(p_{\alpha} : \alpha<\delta)$
is $\leq_{\mathbb P}$-increasing and continuous (i.e. $\alpha<\delta \rightarrow \underset{\beta<\alpha}{\cup} \mathbb{Q}_{p_{1+\beta}} \lessdot \mathbb{Q}_{p_{\alpha}}$),
then $\bar p$ has an upper bound $p_{\delta}$ such that $\bar p \hat{} (p_{\delta})$
is $\leq_{\mathbb P}$-increasing continuous.

c. In clause (b), if $\aleph_2 \leq cf(\delta)$, then $p_{\delta}$
can be chosen as the union of the $p_{\alpha}$s.

d. $\Vdash_{\mathbb P} " \underset{\sim}{ \mathbb Q}$ is ccc and
$\lambda$ is its underlying set of elements$"$.

e. $\Vdash_{\mathbb P} " \Vdash_{\underset{\sim}{\mathbb Q}} "\underset{\sim}{R}$
is a transcendence basis for the reals.

f. Every permutation $g$ of $\lambda$ naturally induces an automorphism
$\hat g$ of $\mathbb P$ and $\underset{\sim}{\mathbb Q}$ which
maps $\underset{\sim}{R}$ to itself.

\textbf{Proof: }a. By clause (b), $\mathbb P$ is $(<\kappa)$-complete,
hence it preserves cardinals and cofinalities $\leq \kappa$. The
rest should be straightforward.

b. As $\underset{\alpha<\delta}{\cup} \mathbb{Q}_{p_{\alpha}}$ is
ccc, it can be extended to a ccc forcing $\mathbb{Q}_{p_{\delta}}$
such that $\underset{\alpha<\delta}{\cup} \mathbb{Q}_{p_{\alpha}} \lessdot \mathbb{Q}_{p_{\delta}}$
and $\Vdash_{\mathbb{Q}_{p_{\delta}}} MA_{\aleph_1}$. As the union
of the $\underset{\sim}{R}_{p_{\alpha}}$ is algebraically independent,
we can extend it to a transcendence basis for the reals.

c. Letting $\mathbb{Q}_{\delta}=\underset{\alpha<\delta}{\cup} \mathbb{Q}_{p_{\alpha}}$,
obviously $\mathbb{Q}_{\delta}$ is ccc. In order to show that $\Vdash_{\mathbb{Q}_{\delta}} MA_{\aleph_1}$,
it's enough to show that for forcing notions of cardinality $\aleph_1$
in $V^{\mathbb{Q}_{\delta}}$. As $\aleph_2 \leq cf(\delta)$, the
names for a given ccc forcing in $V^{\mathbb{Q}_{\delta}}$ and $\aleph_1$-many
of its dense subsets are already $\mathbb{Q}_{\alpha}$-names for
some $\alpha<\delta$, and as $\Vdash_{\mathbb{Q}_{\alpha}} MA_{\aleph_1}$,
we're done. Similarly, every $\mathbb{Q}_{\delta}$-name for a real
is already a $\mathbb{Q}_{\alpha}$-name for some $\alpha<\delta$,
hence $\underset{\alpha<\delta}{\cup} \underset{\sim}{R}_{p_{\alpha}}$
is a $\mathbb{Q}_{\delta}$-name of a transcendence basis.

d. Let $G\subseteq \mathbb P$ be generic over $V$, we shall argue
in $V[G]$. Given $I=\{q_{\alpha} : \alpha<\omega_1\} \subseteq \mathbb Q$,
as $\mathbb P$ is $(<\kappa)$-complete, it doesn't add new sequences
of ordinals of length $\omega_1$, hence $I\in V$. For every $p\in \mathbb P$,
there is some $q\in \mathbb P$ above $p$ such that $I\subseteq \mathbb{Q}_q$.
Therefore, there is some $p\in G$ such that $I\subseteq \mathbb{Q}_p$.
As $\mathbb{Q}_p$ is ccc, there are two elements of $I$ that are
compatible in $\mathbb{Q}_p$ and hence they're compatible in $\mathbb Q$.
It follows that $\mathbb Q$ is ccc. By a similar density argument,
for every $\alpha<\lambda$, there is some $p\in G$ such that $\alpha \in \mathbb{Q}_p$,
hence $\lambda$ is the underlying set of elements of $\mathbb Q$. 

e. As before, we shall argue in $V[G]$ where $G\subseteq \mathbb P$
is generic over $V$. The algebraic independence of $\underset{\sim}{R}$
follows from $G$ being directed. As for the maximality of $\underset{\sim}{R}$,
as before, suppose that $\underset{\sim}{r}$ is a $\mathbb{Q}$-name
for a real, then by a similar argument as in clause (d), there is
$p\in G$ such that $\underset{\sim}{r}$ is a $\mathbb{Q}_p$-name.
As $\underset{\sim}{R}_p$ is a $\mathbb{Q}_p$-name of a transcendence
basis, we're done.

f. This is straightforward. $\square$

\textbf{Definition/Observation 5: }Let $V_1$ be the model $HOD(\mathbb{R}^{<\kappa} \cup \{ \underset{\sim}{R} \} \cup V)$
inside $V^{\mathbb P * \underset{\sim}{Q}}$, then $V_1$ is a model
of $ZF+DC_{< \kappa}$ with the same reals as $V^{\mathbb{P} * \underset{\sim}{Q}}$.
In particular, $V_1$ contains a transcendence basis for the reals
(using Claim 4(e)). $\square$

We shall obtain the desired result by proving that there is no well
ordering of the reals in $V_1$. Before that, we shall prove our main
amalgamation claim:

\textbf{Main amalgamation claim 6: }(A) implies (B) where:

A. a. $\mathbb{Q}_0 \lessdot \mathbb{Q}_l$ $(l=1,2)$.

b. $\Vdash_{\mathbb{Q}_l} "\bar{\underset{\sim}{r}_l}=(\underset{\sim}{r_{l,i}} : i<n_l)$
is algebraically independent over $\mathbb{R}^{V^{\mathbb{Q}_0}}$.

c. $\mathbb Q=\mathbb{Q}_1 \times_{\mathbb{Q}_0} \mathbb{Q}_2$.

B. $\Vdash_{\mathbb Q} "\bar{\underset{\sim}{r_1}} \hat{} \bar{\underset{\sim}{r_2}}$
is algebraically independent over $\mathbb{R}^{V^{\mathbb{Q}_0}}$.

\textbf{Proof: }Assume towards contradiction that there is a counterexample
to the claim. As forcing with $\mathbb Q$ is the same as forcing
with $\mathbb{Q}_0 * ((\mathbb{Q}_1 / \mathbb{Q}_0) \times (\mathbb{Q}_2 / \mathbb{Q}_0))$,
if there is a counterexample to the claim, then by working in $V^{\mathbb{Q}_0}$
we obtain a counterexample where $\mathbb{Q}_0$ is trivial and $\mathbb Q=\mathbb{Q}_1 \times \mathbb{Q}_2$.
Therefore, we may assume wlog that $\mathbb Q=\mathbb{Q}_1 \times \mathbb{Q}_2$
and $\mathbb{Q}_0$ is trivial. We may also assume wlog that it's
forced by $\mathbb Q$ that $\bar{\underset{\sim}{r_1}}$ and $\bar{\underset{\sim}{r_2}}$
form a counterexample (if $(q_1,q_2) \in \mathbb{Q}_1 \times \mathbb{Q}_2$
forces that $\bar{\underset{\sim}{r_1}}$ and $\bar{\underset{\sim}{r_2}}$
form a counterexample, then we can replace $\mathbb{Q}_l$ by $\mathbb{Q}_l \restriction q_l$
for $l=1,2$).

\textbf{Subclaim: }We may assume wlog that $\mathbb{Q}_1$ and $\mathbb{Q}_2$
are Cohen forcing.

Proof of Subclaim: Suppose that $\bar{x}=(\mathbb{Q}_1,\mathbb{Q}_2,\bar{\underset{\sim}{r_1}},\bar{\underset{\sim}{r_2}})$
form a counter example to the amalgamation claim, we shall construct
a counter example $\bar{x'}=(\mathbb{Q}_1',\mathbb{Q}_2',\bar{\underset{\sim}{r_1'}},\bar{\underset{\sim}{r_2'}})$
where $\mathbb{Q}_1',\mathbb{Q}_2'$ are Cohen forcing. As $\bar x$
is a counter example to the claim, there is a nontrivial polynomial
$P=P(x_0,...,x_{n_1-1},y_0,...,y_{n_2-1})$ with coeficients in $\mathbb{R}^V$
and a condition $(p_1,p_2) \in \mathbb{Q}_1 \times \mathbb{Q}_2$
such that $(p_1,p_2) \Vdash_{\mathbb{Q}_1 \times \mathbb{Q}_2} "P(\bar{\underset{\sim}{r_1}},\bar{\underset{\sim}{r_2}})=0"$.
We shall now choose $(\bar{p_{1,n}}, \bar{p_{2,n}}, \bar{a_{1,n}}, \bar{a_{2,n}})$
by induction on $n<\omega$ such that the following conditions hold:

a. $\bar{p_{l,n}}=(p_{l,n,\nu} : \nu \in \omega^n)$ $(l=1,2)$.

b. Each $p_{l,n,\nu}$ is a condition in $\mathbb{Q}_l$ $(l=1,2)$.

c. If $n=m+1$, $l\in \{1,2\}$ and $\nu \in \omega^n$ then $p_{l,m,\nu \restriction m} \leq p_{l,n,\nu}$.

d. $\bar{a_{l,n}}=(a_{l,n,\eta,i}^-,a_{l,n,\eta,i}^+ : \eta \in \omega^n,i<n_l)$.

e. $a_{l,n,\eta,i}^-$ and $a_{l,n,\eta,i}^+$ are rationals such
that $a_{l,n,\eta,i}^+-a_{l,n,\eta,i}^- <\frac{1}{2^n}$.

f. $p_{l,n,\eta} \Vdash_{\mathbb{Q}_l} "\underset{i<n_l}{\wedge}a_{l,n,\eta,i}^- < \underset{\sim}{r_{l,i}} < a_{l,n,\eta,i}^+"$.

g. If $n=m+1$, $\rho \in \omega^m$, $l\in \{1, 2\}$, $((a_i,b_i) : i<n_l)$
is a sequence of pairs of rationals such that $a_i<b_i$ for $i<n_l$
and $p_{l,m,\rho} \nVdash_{\mathbb{Q}_l} "\neg(\underset{i<n_l}{\wedge} a_i<\underset{\sim}{r_{l,i}}<b_i)"$,
then for some $k<\omega$, $p_{l,n,\rho \hat{}(k)} \Vdash_{\mathbb{Q}_l} "\underset{i<n_l}{\wedge} a_i<\underset{\sim}{r_{l,i}}<b_i"$.

h. Moreover, we have $a_i<a_{l,n,\rho \hat{} (k),i}^-<a_{l,n,\rho \hat{} (k),i}^+<b_i$.

i. Moreover, if $n=m+1$ and $\nu_1,\nu_2 \in \omega^m$, then for
some $k_1$ and $k_2$, letting $\rho_l=\nu_l \hat{} (k_l)$ $(l=1,2)$
we have: For all $x_1,...,x_{n_1},y_1,...,y_{n_2}$, if $\underset{i<n_1}{\wedge}a_{l,n,\rho_1,i}^-<x_i<a_{l,n,\rho_1,i}^+$
and $\underset{j<n_2}{\wedge}a_{l,n,\rho_2,j}^-<y_j<a_{l,n,\rho_2,j}^+$
then $-\frac{1}{2^n}<P(x_1,...,x_{n_1-1},y_1,...,y_{n_2-1})<\frac{1}{2^n}$.

j. The $a_{l,n,\eta,i}^-$ are increasing with $\eta$ and the $a_{l,n,\eta,i}^+$
are decreasing with $\eta$.

The induction is straightorward where for clause (i) we use the fact
that $(p_1,p_2) \Vdash_{\mathbb{Q}_1 \times \mathbb{Q}_2} "P(\bar{\underset{\sim}{r_1}},\bar{\underset{\sim}{r_2}})=0"$. 

For $l=1,2$ we define the following objects:

a. $\mathbb{Q}_l'=(\omega^{<\omega}, \leq)$.

b. $\underset{\sim}{\eta_l}$ is the name for the generic real of
$\mathbb{Q}_l'$.

c. For $i<n_l$, $\underset{\sim}{r_{l,i}'}$ is the unique real in
$\underset{n<\omega}{\cap}(a_{l,n,\underset{\sim}{\eta_l \restriction n, i}}^-,a_{l,n,\underset{\sim}{\eta_l \restriction n, i}}^+)$.

Now $\mathbb{Q}_l'$ are equivalent to Cohen forcing, and by clause
(i) of the induction, $\Vdash_{\mathbb{Q}_1' \times \mathbb{Q}_2'} "P(\underset{\sim}{r_1'},\underset{\sim}{r_2'})=0"$.
Therefore, in order to prove the subclaim, it suffices to show that
$\Vdash_{\mathbb{Q}_l'} "\underset{\sim}{r_{l,1}'},...,\underset{\sim}{r_{l,n_1-1}'}$
are algebraically independent over $\mathbb{R}^V"$. Assume towards
contradiction that there is some $\eta \in \mathbb{Q}_l'$ and a nontrivial
polynomial $P_l'(x_0,...,x_{n_l-1})$ such that $\eta \Vdash_{\mathbb{Q}_l'} "P_l'(\underset{\sim}{r_l'})=0"$.
By the assumption on $(\mathbb{Q}_l,\underset{\sim}{r_l})$, letting
$n=lg(\eta)$, $p_{l,n,\eta} \Vdash_{\mathbb{Q}_l} "P_l'(\underset{\sim}{r_l}) \neq 0"$.
Let $G_l \subseteq \mathbb{Q}_l$ be generic over $V$ such that $p_{l,n,\eta} \in G_l$,
so wlog $P_l'(\underset{\sim}{r_l}[G_l])>0$. By continuity, there
are rationals $a_i<b_i$ $(i<n_l)$ such that $V[G]\models "$for
every $x_0,...,x_{n_l-1}$, $\underset{i<n_l}{\wedge}a_i<x_i<b_i \rightarrow P_l'(x_0,...,x_{n_l-1})>0$
and $\underset{\sim}{r_{l,i}}[G_l] \in (a_i,b_i)"$. Therefore, the
first part of the statement holds in $V$ and there is some $q\in G_l$
such that $p_{l,n,\eta} \leq q$ and $q$ forces the second part of
the statement. In particular, $p_{l,n,\eta} \nVdash_{\mathbb{Q}_l} "\neg(\underset{i<n_l}{\wedge}a_i<\underset{\sim}{r_{l,i}}<b_i )"$.
By clause (g) of the induction, there is some $k<\omega$ such that
$p_{l,n+1,\eta \hat{}(k)} \Vdash_{\mathbb{Q}_l} "\underset{i<n_l}{\wedge}a_i<\underset{\sim}{r_{l,i}}<b_i"$
and $a_i<a_{l,n+1,\eta \hat{}(k),i}^-<a_{l,n+1,\eta \hat{}(k),i}^+<b_i$.
Now $\eta \hat{} (k)$ is a condition in $\mathbb{Q}_l'$ that forces
in $\mathbb{Q}_l'$ that $\underset{\sim}{r_{l,i}}' \in (a_i, b_i)$
for all $i<n_l$. It follows that $\eta \hat{}(k)$ forces in $\mathbb{Q}_l'$
that $P_l'(\underset{\sim}{r_{l,0}}',...,\underset{\sim}{r_{l,n_l-1}}')>0$,
contradicting the choise of $\eta$ and $P-l'$. It follows that $\Vdash_{\mathbb{Q}_l'} "\underset{\sim}{r_{l,1}'},...,\underset{\sim}{r_{l,n_1-1}'}$
are algebraically independent over $\mathbb{R}^V"$, which completes
the proof of the subclaim.

We shall now return to the proof of the main amalgamation claim:

Let $\chi \geq \aleph_1$ be large enough and let $N$ be a countable
elementary submodel of $(H(\chi),\in)$ such that $\mathbb{Q}_l,\bar{\underset{\sim}{r_l}} \in N$
$(l=1,2)$. As $\mathbb{Q}_l$ is Cohen, there is a $\mathbb{Q}_l$-name
$\underset{\sim}{\eta_l}$ for a Cohen real over $V$ that generates
the generic for $\mathbb{Q}_l$. For each $l<3$ and $i<n_l$ there
is a Borel function $\bold{B}_{l,i}$ such that $\underset{\sim}{r_{l,i}}=\bold{B}_{l,i}(\underset{\sim}{\eta_l})$,
we may assume that the $\bold{B}_{l,i}$s belong to $N$ as well.
Let $\eta_1' \in V$ be Cohen over $N$, let $G_2 \subseteq \mathbb{Q}_2$
be generic over $V$ and let $\eta_2=\underset{\sim}{\eta_2}[G_2]$.
$\eta_2$ is Cohen over $V$ and is also generic over $N[\eta_1']$.
Therefore, $(\eta_1',\eta_2)$ is generic for $\mathbb{Q}_1 \times \mathbb{Q}_2$
over $N$. As it's forced by $\mathbb{Q}_1 \times \mathbb{Q}_2$ over
$V$ that $\bar{\underset{\sim}{r_1}} \hat{} \bar{\underset{\sim}{r_2}}$
is a counterexample, there is a polynomial $P$ witnessing this, i.e.
$V\models "\Vdash_{\mathbb{Q}_1 \times \mathbb{Q}_2} "P(...,\bold{B}_{1,l}(\underset{\sim}{\eta_1'}),...,...,\bold{B}_{2,l}(\underset{\sim}{\eta_2}),...)=0""$.
By absoluteness, the same stetement holds in $N$. By the genericity
over $N$ of $(\eta_1',\eta_2)$, $N[\eta_1',\eta_2] \models P(...,\bold{B}_{1,l}(\eta_1'),...,...,\bold{B}_{2,l}(\eta_2),...)=0$.
Therefore, there is $p_2 \in G_2 \subseteq \mathbb{Q}_2$ such that
$N[\eta_1'] \models "p_2 \Vdash_{\mathbb{Q}_2} "\bar{\underset{\sim}{r_2}}$
is not algebraically independent over $\mathbb{R}^V$, as witnessed
by $(\bold{B}_{1,l}(\eta_1') : l<n_1)""$, and by absoluteness, the
same holds in $V$. This contradicts assumption (A)(b) and completes
the proof of the claim. $\square$

Before proving the relevant conclusion for $\mathbb P$, we need the
following algebraic observation:

\textbf{Observation 7: }Let $p_1,p_2 \in \mathbb{P}$ and suppose
that $p_1 \leq p_2$. Denote $\mathbb{Q}_{p_l}$ by $\mathbb{Q}_l$
and $\underset{\sim}{R_{p_l}}$ by $\underset{\sim}{R_l}$ $(l=1,2)$.
Then $\Vdash_{\mathbb{Q}_2} "\underset{\sim}{R_2} \setminus \underset{\sim}{R_1}$
is algebraically independent over $\mathbb{R}^{V^{\mathbb{Q}_1}}"$.

\textbf{Proof: }Suppose towards contradiction that there is some $q\in \mathbb{Q}_2$
and $\underset{\sim}{r_0},...,\underset{\sim}{r_{n_2-1}}$ (with no
repetition) such that $q\Vdash_{\mathbb{Q}_2} "\underset{\sim}{r_0},...,\underset{\sim}{r_{n_2-1}} \in \underset{\sim}{R_2} \setminus \underset{\sim}{R_1}$
are not algebraically independent over $\mathbb{R}^{V^{\mathbb{Q}_1}}"$.
By increasing $q$ if necessary, we may assume wlog that there is
a non-trivial polynomial $P(x_0,...,x_{n_2-1})$ over $\mathbb{R}^{V^{\mathbb{Q}_1}}$
such that $q\Vdash_{\mathbb{Q}_2} "P(\underset{\sim}{r_0},...,\underset{\sim}{r_{n_2-1}})=0"$.
Therefore, there are $\mathbb{Q}_1$-names of reals $\underset{\sim}{s_0},...,\underset{\sim}{s_{n_1-1}}$
and a polynomial $Q(x_0,...,x_{n_2-1},y_0,...,y_{n_1-1})$ over the
rationals such that $q\Vdash_{\mathbb{Q}_2} "Q(x_0,...,x_{n_2-1}, \underset{\sim}{s_0},...,\underset{\sim}{s_{n_1-1}})=P(x_0,...,x_{n_2-1})"$.
Recalling that $\underset{\sim}{R_1}$ is a $\mathbb{Q}_1$-name of
a transcendence basis over the rationals, then by increasing $q$
if necessary, there are $\mathbb{Q}_1$-names of reals $\underset{\sim}{t_0},...,\underset{\sim}{t_{n_0-1}}$
such that $q\Vdash_{\mathbb{Q}_2} "\underset{\sim}{t_0},...,\underset{\sim}{t_{n_0-1}} \in \underset{\sim}{R_1}$
(with no repetition)$"$ and $q\Vdash_{\mathbb{Q}_2} "\underset{\sim}{s_0},...,\underset{\sim}{s_{n_1-1}}$
are algebraic over $\mathbb{Q}[\underset{\sim}{t_0},...,\underset{\sim}{t_{n_0-1}}]"$
(here $\mathbb Q$ denotes the field of rational numbers). It follows
that $q\Vdash_{\mathbb{Q}_2} "\{\underset{\sim}{t_0},...,\underset{\sim}{t_{n_0-1}},\underset{\sim}{r_0},...\underset{\sim}{r_{n_2-1}} \} \subseteq \underset{\sim}{R_2}$
is not algebraically independent over the rationals$"$. By the choice
of the $\underset{\sim}{t_i}$s and the $\underset{\sim}{r_i}$s,
$q\Vdash_{\mathbb{Q}_2} "\underset{\sim}{t_0},...,\underset{\sim}{t_{n_0-1}},\underset{\sim}{r_0},...,\underset{\sim}{r_{n_2-1}}$
are without repetition''. Together, we get a contradiction to the
definition of the conditions in $\mathbb{P}$ and the fact that $p_2 \in \mathbb P$.
$\square$

\textbf{Conclusion 8: }Suppose that $p_1, p_2 \in \mathbb{P}$ such
that $p_1 \leq p_2$. Let $g$ be a permutation of $\lambda$ of order
$2$ such that $g\restriction u_{p_1}=id$ and $g''(u_{p_2}) \cap u_{p_2}=u_{p_1}$,
and let $p_3=\hat{g}(p_2)$. Then there is $q\in \mathbb P$ such
that $p_2,p_3 \leq q$ and $\mathbb{Q}_{p_2} \times_{\mathbb{Q}_{p_1}} \mathbb{Q}_{p_3} \lessdot \mathbb{Q}_q$.

\textbf{Proof: }Let $\mathbb{Q}= \mathbb{Q}_{p_2} \times_{\mathbb{Q}_{p_1}} \mathbb{Q}_{p_3}$.
As $\mathbb{Q}_1$ is ccc and $\Vdash_{\mathbb{Q}_1} "MA_{\aleph_1}+\mathbb{Q}_2 / \mathbb{Q}_1  \models ccc + \mathbb{Q}_3 / \mathbb{Q}_1 \models ccc"$,
it follows that $\mathbb Q$ is ccc (see e.g. {[}HwSh:1090{]} for
details). By the previous observation, for $l=2,3$, $\Vdash_{\mathbb{Q}_{p_l}} "\underset{\sim}{R_{p_l}} \setminus \underset{\sim}{R_{p_1}}$
is algebraically independent over $\mathbb{R}^{V^{\mathbb{Q}_{p_1}}}"$.
Therefore, by Claim 6, $\Vdash_{\mathbb Q} "(\underset{\sim}{R_{p_2}} \setminus \underset{\sim}{R_{p_1}}) \cup (\underset{\sim}{R_{p_3}} \setminus \underset{\sim}{R_{p_1}})$
is algebraically independent over $\mathbb{R}^{V^{\mathbb{Q}_1}}"$.
It follows that $\Vdash_{\mathbb{Q}} "\underset{\sim}{R_{p_2}} \cup \underset{\sim}{R_{p_3}}=\underset{\sim}{R_{p_1}} \cup (\underset{\sim}{R_{p_2}} \setminus \underset{\sim}{R_{p_1}}) \cup (\underset{\sim}{R_{p_3}} \setminus \underset{\sim}{R_{p_1}})$
is algebraically independent over the rationals$"$ (recall that if
$\{ \alpha_0,...,\alpha_{n-1}\}$ are algebraically independent over
the rationals and $\{ \alpha_n,...,\alpha_{m-1}\}$ are algebraically
independent over a field $\mathbb F$ containing $\mathbb Q \cup \{ \alpha_0,...,\alpha_{n-1}\}$,
then $\{ \alpha_0,...,\alpha_{m-1} \}$ are algebraically independent
over the rationals). By Hypothesis 1, there is a ccc forcing $\mathbb{Q}_{q}$
such that $\mathbb Q \lessdot \mathbb{Q}_{q}$, $\Vdash_{\mathbb{Q}_{q}} MA_{\aleph_1}$
and $|\mathbb{Q}_{q}|=u_{q}$ for some $u_{q} \in [\lambda]^{<\kappa}$.
As $\Vdash_{\mathbb{Q}_q} "\underset{\sim}{R_{p_2}} \cup \underset{\sim}{R_{p_3}}$
are algebraically independent over the rationals$"$, there is a set
$\underset{\sim}{R_q}$ of $\mathbb{Q}_q$-names of reals such that
$\underset{\sim}{R_{p_2}} \cup \underset{\sim}{R_{p_3}} \subseteq \underset{\sim}{R_{q}}$
and $\Vdash_{\mathbb{Q}_q} " \underset{\sim}{R_q}$ is a transcendence
basis for the reals$"$. Now let $q=(u_q,\mathbb{Q}_q,\underset{\sim}{R_q})$,
it's easy to verify that $q$ is as required. $\square$

Recalling Observation 5, we shall complete the proof of the main result
of the paper by proving the following claim:

\textbf{Claim 9: }There is no well-ordering of the reals in $V_1$.

\textbf{Proof: }Assume towards contradiction that there are $(p_1,r_1)\in \mathbb P * \underset{\sim}{\mathbb Q}$
such that, over $V$, $(p_1,r_1) \Vdash_{\mathbb{P} * \underset{\sim}{\mathbb Q}} "\underset{\sim}{f}$
is a one-to-one function from $\mathbb R$ to $Ord"$ and such that
$\underset{\sim}{f}$ is definable from $\underset{\sim}{R}$ and
a sequence $(\underset{\sim}{\eta_{\epsilon}} : \epsilon<\epsilon(*))$
where $\epsilon(*)<\kappa$ and wlog each $\underset{\sim}{\eta}_{\epsilon}$
is a $\mathbb{Q}_{p_1}$ name for a real (by a similar argument as
in claims 4(d) and 4(e), we can always extend $p_1$ to make this
true). Choose $(p_2,r_2) \geq (p_1, r_1)$ and a name of a real $\underset{\sim}{r}$
such that $(p_2,r_2) \Vdash_{\mathbb{P} * \underset{\sim}{\mathbb Q}} "\underset{\sim}{r} \in \mathbb{R}^{V^{\mathbb{Q}_{p_2}}} \setminus \mathbb{R}^{V^{\mathbb{Q}_{p_1}}}"$,
wlog $r_2 \in \mathbb{Q}_{p_2}$, and by extending the condition if
necessary, we may assume wlog that $(p_2,r_2)$ forces a value $\gamma$
to $\underset{\sim}{f}(\underset{\sim}{r})$.

Let $g$ be a permutation of $\lambda$ of order $2$ such that $g\restriction u_{p_1}=id$
and $g''(u_{p_2}) \cap u_{p_2}=u_{p_1}$. We shall denote both of
the induced automorphisms on $\mathbb P$ and $\mathbb Q$ by $\hat{g}$.
Clearly, $\hat{g}(p_1)=p_1$. Let $p_3=\hat{g}(p_2)$ and $r_3=\hat{g}(r_2)$.
By the previous claims, there is $q\in \mathbb P$ such that $p_2,p_3 \leq q$
and $\mathbb{Q}_{p_2} \times_{\mathbb{Q}_{p_1}} \mathbb{Q}_{p_3} \lessdot \mathbb{Q}_q$,
and by the construction of the amalgamation, there is $r\in \mathbb{Q}_q$
above $r_2$ and $r_3$. As $\Vdash_{\mathbb{P} * \underset{\sim}{\mathbb Q}} "\mathbb{R}^{V^{\mathbb{Q}_{p_2}}} \cap \mathbb{R}^{V^{\mathbb{Q}_{p_3}}}=\mathbb{R}^{V^{\mathbb{Q}_{p_1}}}"$,
it follows that $(q,r) \Vdash_{\mathbb{P} * \underset{\sim}{\mathbb Q}} "\underset{\sim}{r} \neq g(\underset{\sim}{r})"$.
As $(p_2,r_2) \leq (q,r)$, $(q,r) \Vdash_{\mathbb{P} * \underset{\sim}{\mathbb Q}} "\underset{\sim}{f}(\underset{\sim}{r})=\gamma"$.
Recalling that $\underset{\sim}{f}$ is forced to be injective, we
shall arrive at a contradiction by showing that $(q,r) \Vdash_{\mathbb{P} * \underset{\sim}{Q}} "\underset{\sim}{f}(\hat{g}(\underset{\sim}{r}))=\gamma"$.
It's enough to show that the statement is forced by $(p_3,r_3)=(\hat{g}(p_2),\hat{g}(r_2))$,
and in order to show that, it suffices to show that $\underset{\sim}{f}=\hat{g}(\underset{\sim}{f})$.
Recalling that each $\underset{\sim}{\eta_{\epsilon}}$ in the definition
of $\underset{\sim}{f}$ is a $\mathbb{Q}_{p_1}$-name and that $g$
is the identity on $u_{p_1}$, it follows that $\hat{g}(\underset{\sim}{\eta_{\epsilon}})=\underset{\sim}{\eta_{\epsilon}}$.
By Claim 4(f), $\underset{\sim}{R}$ is preserved by $\hat{g}$. As
$\underset{\sim}{f}$ is definable from $\underset{\sim}{R}$ and
$(\underset{\sim}{\eta_{\epsilon}} : \epsilon< \epsilon(*))$, it
follows that $\hat{g}(\underset{\sim}{f})=\underset{\sim}{f}$. This
completes the proof of the claim. $\square$

\textbf{\large References}{\large \par}

{[}HwSh:1090{]} Haim Horowitz and Saharon Shelah, Can you take Toernquist's
inaccessible away?, arXiv:1605.02419

{[}HwSh:1113{]} Haim Horowitz and Saharon Shelah, Madness and regularity
properties, arXiv:1704.08327

{[}LaZa1{]} Paul Larson and Jindrich Zapletal, Canonical models for
fragments of the axiom of choice, Journal of Symbolic Logic, 82:489-509,
2017

{[}LaZa2{]} Paul Larson and Jindrich Zapletal, Geometric set theory,
preprint

{[}Sh:218{]} Saharon Shelah, On measure and category, Israel J. Math.
52 (1985) 110-114

$\\$

(Haim Horowitz) Department of Mathematics

University of Toronto

Bahen Centre, 40 St. George St., Room 6290

Toronto, Ontario, Canada M5S 2E4

E-mail address: haim@math.toronto.edu

$\\$

(Saharon Shelah) Einstein Institute of Mathematics

Edmond J. Safra Campus,

The Hebrew University of Jerusalem.

Givat Ram, Jerusalem 91904, Israel.

Department of Mathematics

Hill Center - Busch Campus,

Rutgers, The State University of New Jersey.

110 Frelinghuysen Road, Piscataway, NJ 08854-8019 USA

E-mail address: shelah@math.huji.ac.il
\end{document}